\documentclass[11pt]{article}

\def\mathbb{\Bbb}
\usepackage{amsfonts,latexsym,amstext,theorem,euscript,array,enumerate,mathrsfs}
\usepackage{amsmath}
\usepackage{amssymb}

\usepackage{comment}

\newtheorem{Theorem}{Theorem}[part]

\newtheorem{Lemma}{Lemma}[part]
\newtheorem{Corollary}{Corollary}[part]



\def\qed{{\hfill\hbox{\enspace${ \square}$}} \smallskip}
\def\sqr#1#2{{\vcenter{\vbox{\hrule height .#2pt \hbox{\vrule
 width .#2pt height#1pt \kern#1pt \vrule
width .#2pt} \hrule height .#2pt}}}}
\def\square{\mathchoice\sqr54\sqr54\sqr{4.1}3\sqr{3.5}3}

\def\ds{\begin{displaystyle}}
\def\eds{\end{displaystyle}}

\def\<{\langle }
\def\>{\rangle }

\def\R{\mathbb R}

\def\E{\mathbb E}
\def\P{\mathbb P}

\def\F{\mathbb F}

\def\cala{{\cal A}}

\def\calf{{\cal F}}

\def\calh{{\cal H}}

\def\calt{{\cal T}}

\def\calv{{\cal V}}
\def\call{{\cal L}}
\def\cals{{\cal S}}

\def\beqs{\begin{eqnarray*}}
\def\enqs{\end{eqnarray*}}
\def\beq{\begin{eqnarray}}
\def\enq{\end{eqnarray}}

\title{Non-Markovian optimal stopping problems \\and
constrained BSDEs with jump
}

\date{}

\author{Marco Fuhrman
\\Politecnico di Milano,
Dipartimento di Matematica\\
via Bonardi 9, 20133 Milano, Italy\\
marco.fuhrman@polimi.it
\\
\\
Huy\^{e}n Pham
\\
LPMA - Universit\'{e} Paris Diderot \\
Batiment Sophie Germain, Case 7012 \\
13 rue Albert Einstein, 75205 Paris Cedex 13\\
and CREST-ENSAE \\
pham@math.univ-paris-diderot.fr
\\
\\
Federica Zeni
\\Politecnico di Milano,
Dipartimento di Matematica\\
via Bonardi 9, 20133 Milano, Italy\\
federica.zeni@mail.polimi.it
}

\begin{document}

 \maketitle

 \begin{abstract}
 We consider a non-Markovian optimal stopping problem on finite horizon.
 We prove that the value process can be represented by means of
 a backward stochastic differential equation (BSDE), defined on an enlarged
 probability space, containing a stochastic integral
 having  a
 one-jump  point process as integrator
 and an (unknown) process with a sign constraint as  integrand.
This provides an alternative representation with respect to the classical
one given by a reflected BSDE. The connection between the two BSDEs
is also clarified. Finally, we prove that the value of the optimal
stopping problem is the same as the value of an auxiliary
optimization  problem where the intensity of
the point process is controlled.
 \end{abstract}


\noindent {\bf MSC Classification (2010):} 60H10, 60G40, 93E20.


\section{Introduction}
Let $(\Omega,\calf,\P)$ be a complete probability space and let $\F=(\calf_t)_{t\ge 0}$
be the natural augmented filtration generated by an $m$-dimensional
standard Brownian motion $W$. For given $T>0$ we denote
$L^2_T=L^2(\Omega,\calf_T,\P)$ and introduce the following spaces
of processes.
\begin{enumerate}
\item $\calh^2=\{Z:\Omega\times [0,T]\to \R^m,\,\F$-predictable,
$\|Z\|^2_{\calh^2}=\E\int_0^T|Z_s|^2ds<\infty\}$;
\item
$\cals^2=\{Y:\Omega\times [0,T]\to \R ,\,\F$-adapted and c\`adl\`ag,
$\|Y\|^2_{\cals^2}=\E\sup_{t\in [0,T]}|Y_s|^2<\infty\}$;
\item
$\cala^2=\{K\in \cals^2,\,\F$-predictable, nondecreasing, $K_0=0\}$;
\item
$\cals_c^2=\{Y\in  \cals^2 $ with continuous paths$\}$;
\item
$\cala_c^2=\{K\in  \cala^2 $ with continuous paths$\}$.
\end{enumerate}
We suppose we are given
\begin{equation}\label{ipotesiriflessa}
f\in\calh^2,\quad h\in \cals_c^2,\quad
\xi\in L^2_T,\quad \text{satisfying}\quad\xi\ge h_T.
\end{equation}
We wish to characterize the process defined, for every $t\in[0,T]$, by
$$
I_t=\mathop{\rm ess\,sup}_{\tau\in\calt_t(\F)}
\E\left[ \int_t^{T\wedge \tau}f_s\,ds +
h_\tau\,1_{\tau< T}+
\xi\,1_{\tau\ge T}\,\bigg|\, \calf_t\right],
$$
where $\calt_t(\F)$ denotes the set of $\F$-stopping times $\tau\ge t$.
Thus, $I$ is the value process of a non-Markovian
optimal stopping problem with cost functions $f,h,\xi$.
In \cite{ElKetal97} the process $I$ is
described by means of an associated reflected backward stochastic differential
equation (BSDE), namely it is proved that
there exists a unique $(Y,Z,K)\in \cals_c^2\times \calh^2\times
\cala_c^2$ such that, $\P$-a.s.
\begin{eqnarray}\label{BSDEriflessa}
Y_t+\int_t^TZ_s\,dW_s&=&
\xi +\int_t^Tf_s\,ds+K_T-K_s,
\\
Y_t\ge h_t, &&
\int_0^T(Y_s-h_s)\,dK_s=0, \qquad t\in [0,T],
\label{BSDEriflessavincolo}
\end{eqnarray}
and that, for every  $t\in [0,T]$, we have $I_t=Y_t$ $\P$-a.s.

It is our purpose to present another representation of the process
$I$ by means of a different BSDE, defined on an enlarged
probability space, containing a jump part  and
involving sign constraints.
Besides its intrinsic interest, this result may lead to new methods
for the numerical approximation of the value process, based
on numerical schemes designed
to approximate the solution to the modified BSDE. In the context of
a classical Markovian optimal stopping problem, this
may give rise to new computational methods for the
corresponding variational inequality as studied in \cite{BeLi}.

We use a  randomization method, which
consists in replacing the stopping time $\tau$ by a random variable $\eta$
independent of the Brownian motion
and in formulating an auxiliary
optimization problem where we can control the intensity of
the (single jump) point process $N_t=1_{\eta\le t}$.
The auxiliary randomized problem turns  out to have the same
value process as the original one. This approach is in the same spirit as in
\cite{khaetal10},  \cite{khph},
\cite{ElKh},
\cite{ElKh-b},    \cite{FuPh} where
BSDEs with barriers and
optimization
problems with switching, impulse control and continuous control were
considered.

\section{Statement of the main results}

We are given $(\Omega,\calf,\P)$,  $\F=(\calf_t)_{t\ge 0}$,
 $W$, $T$ as before, as well as $f,h,\xi$
satisfying
\eqref{ipotesiriflessa}.
Let $\eta$ be an exponentially distributed random variable with unit mean,
defined in another probability space  $(\Omega',\calf',\P')$.
 Define $\bar\Omega= \Omega\times \Omega'$
 and let $(\bar\Omega, \bar\calf,\bar\P)$ be the completion of
 $(\bar\Omega, \calf\otimes\calf',\P \otimes \P')$.
 All the random elements $W,f,h,\xi,\eta$ have natural extensions to $\bar\Omega$,
 denoted by the same symbols.  Define
 $$
 N_t=1_{\eta\le t}, \qquad A_t=t\wedge \eta,
 $$
 and let  $\bar\F=(\bar\calf_t)_{t\ge 0}$
be the $\bar\P$-augmented filtration generated by   $(W,N)$.
Under $\bar\P$,
  $A$ is the $\bar\F$-compensator (i.e., the dual predictable
 projection) of $N$, $W$ is an $\bar\F$-Brownian motion
 independent of $N$
and \eqref{ipotesiriflessa} still holds provided
$\calh^2$,  $\cals_c^2$,
$L^2_T$ (as well as $\cala^2$ etc.) are understood
with respect to $(\bar\Omega,\bar\calf, \bar\P)$ and $\bar\F$
as we will do. We also define
$$\call^2=\{U:\bar\Omega\times [0,T]\to \R,\,\bar\F{\rm -predictable},
\quad
\|U\|^2_{\call^2}=\bar\E\int_0^T|U_s|^2dA_s=
\bar\E\int_0^T|U_s|^2dN_s<\infty\}.
$$
We will consider the BSDE
\begin{equation}\label{BSDEconstrained}
\bar Y_t+\int_t^T\bar Z_s\,dW_s+\int_{(t,T]}\bar U_s\,dN_s=
\xi\,1_{\eta\ge T} +\int_t^Tf_s\,1_{[0,\eta]}(s)\,ds+
\int_{(t,T]}h_s\,dN_s+
\bar K_T-\bar K_t,
 \quad t\in [0,T],
\end{equation}
with the constraint
\begin{equation}\label{constraint}
U_t\le 0, \qquad dA_t(\bar\omega)\,\bar\P(d\bar\omega)-a.s.
\end{equation}
We say that
 a quadruple  $(\bar Y, \bar Z, \bar U, \bar K)$
 is  a solution to this BSDE if it belongs to $\cals^2\times \calh^2\times
\call^2\times \cala^2$,  \eqref{BSDEconstrained} holds      $\bar\P$-a.s.,
and \eqref{constraint} is satisfied.
We say that
$(\bar Y, \bar Z, \bar U, \bar K)$ is minimal if for any other solution
$(\bar Y', \bar Z', \bar U', \bar K')$ we have, $\bar\P$-a.s,
 $\bar Y_t\le \bar Y_t'$ for all $t\in [0,T]$.

 Our first main result shows the existence of a minimal solution to the BSDE
 with sign constraint and makes the connection with
reflected BSDEs.

\begin{Theorem}\label{BSDEwellposed}
Under \eqref{ipotesiriflessa} there exists a unique minimal
solution $(\bar Y, \bar Z, \bar U, \bar K)$
to \eqref{BSDEconstrained}-\eqref{constraint}. It
can be defined  starting from the solution
$(Y,Z,K)$
to the reflected BSDE
\eqref{BSDEriflessa}-\eqref{BSDEriflessavincolo} and setting,
  for $\bar \omega=(\omega,\omega')$,
$t\in [0,T]$,
\begin{eqnarray}
   \bar Y_t(\bar\omega)=   Y_t( \omega)1_{t<  \eta(\omega')}, &&
    \bar Z_t(\bar\omega)= Z_t( \omega)1_{t\le   \eta(\omega')},
    \label{identifuno}\\
   \bar U_t(\bar\omega)= (h_t(\omega)-Y_{t}(\omega))
 1_{t\le   \eta(\omega')},
 &&
  \bar K_t(\bar\omega)= K_{t \wedge \eta(\omega')}(\omega).
  \label{identifdue}
\end{eqnarray}
\end{Theorem}

Now we formulate an auxiliary optimization problem.
 Let $\calv=\{\nu:\bar\Omega\times [0,\infty)\to (0,\infty),$
$\bar\F$-predictable and bounded$\}$. For $\nu\in\calv$ define
$$
L_t^\nu=\exp\left(
\int_0^t (1-\nu_s)\,dA_s+\int_0^t\log\nu_s\,dN_s\right)
= \exp\left(
\int_0^{t\wedge \eta}(1-\nu_s)\,ds\right)(1_{t<\eta}+\nu_\eta1_{t\ge\eta}).
$$
Since $\nu$ is bounded, $L^\nu$ is an $\bar\F$-martingale
on $[0,T]$ under $\bar\P$ and
we can define an equivalent
probability $\bar \P_\nu$ on $(\bar\Omega, \bar \calf)$ setting
$\bar \P_\nu(d\bar \omega)=L_t^\nu(\bar\omega)\,\bar\P(d\bar\omega)$.
By a theorem of  Girsanov type
(Theorem 4.5 in \cite{ja}) on $[0,T]$ the $\bar\F$-compensator of $N$ under $\bar\P_\nu$
is $\int_0^t \nu_s\,dA_s$, $t\in [0,T]$, and $W$ remains a Brownian motion
under $\bar\P_\nu$. We wish to characterize the value process $J$
defined, for every $t\in[0,T]$, by
\begin{equation}\label{dualrepresentation}
    J_t=\mathop{\rm ess\,sup}_{\nu\in\calv}
\bar \E_\nu\left[ \int_{t\wedge \eta}^{T\wedge \eta}f_s\,ds +
h_\eta\,1_{t<\eta< T}+
\xi\,1_{\eta\ge T}\,\bigg|\, \bar\calf_t\right].
\end{equation}

Our second result provides a dual representation in terms
of control intensity of the minimal solution to the BSDE with sign constraint.

\begin{Theorem}\label{dualproblem} Under \eqref{ipotesiriflessa}, let
 $(\bar Y, \bar Z, \bar U, \bar K)$
be
the  minimal solution
to \eqref{BSDEconstrained}-\eqref{constraint}.
Then, for every $t\in[0,T]$, we have $\bar Y_t=J_t$ $\bar\P$-a.s.
\end{Theorem}

The equalities $J_0=\bar Y_0=Y_0 =I_0$ immediately give  the following
corollary.
\begin{Corollary}
Under \eqref{ipotesiriflessa}, let
 $(\bar Y, \bar Z, \bar U, \bar K)$
be
the  minimal solution
to \eqref{BSDEconstrained}-\eqref{constraint}. Then
$$
\bar Y_0=
     \sup_{\tau\in\calt_0(\F)}
\E\left[ \int_0^{T\wedge \tau}f_s\,ds +
h_\tau\,1_{\tau< T}+
\xi\,1_{\tau\ge T}\right]
    = \sup_{\nu\in\calv}
\bar \E_\nu\left[ \int_0^{T\wedge \eta}f_s\,ds +
h_\eta\,1_{\eta< T}+
\xi\,1_{\eta\ge T}\right].
$$
\end{Corollary}

\section{Proofs}

\noindent {\bf Proof of Theorem \ref{BSDEwellposed}.}
Uniqueness of the minimal solution is not difficult and it is established as in \cite{khph},
Remark 2.1.

Let $(Y,Z,K)\in \cals_c^2\times \calh^2\times
\cala_c^2$ be the solution to
 \eqref{BSDEriflessa}-\eqref{BSDEriflessavincolo}, and let
 $(\bar Y, \bar Z, \bar U, \bar K)$ be defined by
\eqref{identifuno},
  \eqref{identifdue}. Clearly
it belongs to $\cals^2\times \calh^2\times
\call^2\times \cala^2$ and
the constraint  \eqref{constraint} is satisfied
 due to the reflection inequality in \eqref{BSDEriflessavincolo}.
The fact that it satisfies equation
\eqref{BSDEconstrained} can be proved by direct substitution,
by considering the three disjoint events
$\{\eta >T\}$, $\{0\le t<\eta <T \}$, $\{0<\eta <T, \eta\le t\le T\}$,
whose union is $\bar\Omega$, $\bar\P$-a.s.

Indeed, on $\{\eta >T\}$ we have $Z_s=\bar Z_s$ for every $s\in [0,T]$
and, by the local property
of the stochastic integral,  $\int_t^T\bar Z_s\,dW_s=\int_t^T Z_s\,dW_s$,
$\bar\P$-a.s. and \eqref{BSDEconstrained} reduces to
\eqref{BSDEriflessa}.

On $\{0\le t<\eta <T \}$ \eqref{BSDEconstrained} reduces to
$$
\bar Y_t+\int_t^T\bar Z_s\,dW_s+ \bar U_\eta=
 \int_t^\eta f_s\, ds+
 h_\eta +
\bar K_T-\bar K_t, \qquad \bar\P-a.s.;
$$
since $\int_t^T\bar Z_s\,dW_s=\int_t^\eta   Z_s\,dW_s$ $\P$-a.s.,
$h_\eta- \bar U_\eta=Y_\eta$
 and,
on the set $\{0\le t<\eta <T \}$, $\bar Y_t=Y_t$ and
$\bar K_T-\bar K_t= K_\eta-  K_t$, this reduces to
$$
 Y_t+\int_t^\eta  Z_s\,dW_s =
 \int_t^\eta f_s\, ds+
 Y_\eta +
  K_\eta-  K_t, \qquad \bar \P-a.s.
$$
which again holds by \eqref{BSDEriflessa}.

Finally,
on $\{0<\eta <T, \eta\le t\le T\}$ the verification of
\eqref{BSDEconstrained}
is trivial, so we have proved that
$(\bar Y, \bar Z, \bar U, \bar K)$ is indeed a solution.

Its minimality property   will be proved later.
\qed

To proceed further we recall a result from \cite{ElKetal97}:
for every integer $n\ge 1$, let
  $(Y^n,Z^n)\in \cals_c^2\times \calh^2$ denote the unique solution to
  the penalized BSDE
\begin{eqnarray}\label{BSDEriflessapenalizzata}
Y_t^n+\int_t^TZ_s^n\,dW_s&=&
\xi +\int_t^Tf_s\,ds+n
\int_t^T(Y_s^n-h_s)^-\,ds , \qquad t\in [0,T];
\end{eqnarray}
then, setting $K_t^n=n
\int_0^t(Y_s^n-h_s)^-\,ds$, the triple $(Y^n,Z^n,K^n)$
converges   in $\cals_c^2\times \calh^2\times
\cala_c^2$ to the solution $(Y,Z,K)$ to  \eqref{BSDEriflessa}-\eqref{BSDEriflessavincolo}.

Define
$$
   \bar Y_t^n(\bar\omega)=   Y_t^n( \omega)1_{t<  \eta(\omega')},
   \quad
    \bar Z_t^n(\bar\omega)= Z_t^n( \omega)1_{t\le  \eta(\omega')},
    \quad
   \bar U_t^n(\bar\omega)= (h_t(\omega)-Y_{t}^n(\omega))
 1_{t\le   \eta(\omega')},
$$
and note that $
   \bar Y^n\to
   \bar Y$ in $\cals^2$.

\begin{Lemma}
$(\bar Y^n, \bar Z^n, \bar U^n)$ is the
 unique solution
in $\cals^2\times \calh^2\times \call^2$ to the BSDE:
$\bar\P$-a.s.,
\begin{eqnarray}\label{BSDEconstrainedpenalizzata}
\bar Y_t^n+\int_t^T\bar Z_s^n\,dW_s+\int_{(t,T]}\bar U_s^n\,dN_s
&=&
\xi\,1_{\eta\ge T} +\int_t^Tf_s\,1_{[0,\eta]}(s)\,ds
\\
&&+
\int_{(t,T]}h_s\,dN_s+
n
\int_t^T(\bar U_s^n)^+1_{[0,\eta]}(s) \,ds , \qquad t\in [0,T].
\nonumber
\end{eqnarray}
\end{Lemma}

\noindent {\bf Proof.}
$(\bar Y^n, \bar Z^n, \bar U^n)$
belongs to $\cals^2\times \calh^2\times \call^2$ and,
proceeding as in the
proof of Theorem \ref{BSDEwellposed} above, one verifies
by direct substitution that
\eqref{BSDEconstrainedpenalizzata} holds,
as a consequence of equation
\eqref{BSDEriflessapenalizzata}.
The uniqueness (which is not needed in the sequel)
follows from
the results in
\cite{Bech}.  \qed

We will identify $\bar Y^n$ with the value process of a
penalized optimization problem.
Let $\calv_n$ denote the set
of all $\nu\in\calv$ taking values in $ (0,n]$
and let us define
(compare  with \eqref{dualrepresentation})
\begin{equation}\label{dualrepresentationpenalized}
    J^n_t=\mathop{\rm ess\,sup}_{\nu\in\calv_n}
\bar \E_\nu\left[ \int_{t\wedge \eta}^{T\wedge \eta}f_s\,ds +
h_\eta\,1_{t<\eta< T}+
\xi\,1_{\eta\ge T}\,\bigg|\, \bar\calf_t\right].
\end{equation}

\begin{Lemma}
For every $t\in[0,T]$, we have $\bar Y_t^n=J_t^n$ $\bar\P$-a.s.
\end{Lemma}

\noindent {\bf Proof.}
We fix any $\nu\in\calv_n$ and
recall that, under the probability
$\bar\P_\nu$,  $W$ is a Brownian motion and the compensator
of $N$
on $[0,T]$
is $\int_0^t \nu_s\,dA_s$, $t\in [0,T]$.
Taking the conditional expectation given $\bar\calf_t$
in \eqref{BSDEconstrainedpenalizzata}
we obtain
\begin{eqnarray}\nonumber
\bar Y_t^n+\bar \E_\nu\left[\int_{(t,T]}\bar U_s^n\,\nu_s\,dA_s
\,\bigg|\, \bar\calf_t\right]
&=&\bar \E_\nu\left[
\xi\,1_{\eta\ge T} +\int_t^Tf_s\,1_{[0,\eta]}(s)\,ds
+
\int_{(t,T]}h_s\,dN_s\,\bigg|\, \bar\calf_t\right]
\\
&&+
\bar \E_\nu\left[
n
\int_t^T(\bar U_s^n)^+1_{[0,\eta]}(s) \,ds \,\bigg|\, \bar\calf_t\right].
\nonumber
\end{eqnarray}
We note that
$\int_{(t,T]}h_s\,dN_s=h_\eta\,1_{t<\eta\le T}=h_\eta\,1_{t<\eta< T}$
$\bar\P_\nu$-a.s., since
  $\eta\ne T$ $\bar\P$-a.s. and hence
$\bar\P_\nu$-a.s.
Since $dA_s=1_{[0,\eta]}(s) \,ds$ we have
\begin{equation}\label{relfond}
    \bar Y_t^n=
\bar \E_\nu\left[ \xi\,1_{\eta\ge T}+\int_{t\wedge \eta}^{T\wedge \eta}f_s\,ds +
h_\eta\,1_{t<\eta< T}\,\bigg|\, \bar\calf_t\right]
+
\bar \E_\nu\left[
\int_t^T(n(\bar U_s^n)^+
-\bar U_s^n\,\nu_s
)1_{[0,\eta]}(s) \,ds \,\bigg|\, \bar\calf_t\right].
\end{equation}
Since $nU^+
-U\,\nu\ge 0$ for every real number $U$ and every $\nu\in (0,n]$ we obtain
$$
\bar Y_t^n\ge
\bar \E_\nu\left[ \xi\,1_{\eta\ge T}+\int_{t\wedge \eta}^{T\wedge \eta}f_s\,ds +
h_\eta\,1_{t<\eta< T}\,\bigg|\, \bar\calf_t\right]
$$
for arbitrary $\nu\in\calv_n$, which implies
$\bar Y_t^n\ge J_t^n$.  On the other hand, setting $\nu^\epsilon_s
= n\,1_{\bar U_s^n>0}+\epsilon \,1_{-1\le \bar U_s^n\le 0}  -
\epsilon \,(\bar U_s^n)^{-1}\,1_{\bar U_s^n<-1}$,
we have $\nu^\epsilon\in\calv_n$ for $0<\epsilon\le 1$ and
$n(\bar U_s^n)^+
-\bar U_s^n\,\nu_s\le \epsilon$. Choosing $\nu=\nu^\epsilon$
in \eqref{relfond} we obtain
$$
    \bar Y_t^n\le
\bar \E_{\nu^\epsilon}\left[ \xi\,1_{\eta\ge T}
+\int_{t\wedge \eta}^{T\wedge \eta}f_s\,ds +
h_\eta\,1_{t<\eta< T}\,\bigg|\, \bar\calf_t\right]
+\epsilon \,T
\le J_t^n +\epsilon \,T
$$
and we have the desired conclusion.
\qed

\noindent {\bf Proof of Theorem \ref{dualproblem}.}
Let
$(\bar Y', \bar Z', \bar U', \bar K')$
be any (not necessarily minimal) solution to \eqref{BSDEconstrained}-\eqref{constraint}.
Since $\bar U'$ is nonpositive and $\bar K'$ is nondecreasing we have
$$
\bar Y_t'+\int_t^T\bar Z_s'\,dW_s\ge
\xi\,1_{\eta\ge T} +\int_t^Tf_s\,1_{[0,\eta]}(s)\,ds+
\int_{(t,T]}h_s\,dN_s=\xi\,1_{\eta\ge T}+
\int_{t\wedge \eta}^{T\wedge \eta}f_s\,ds +
h_\eta\,1_{t<\eta\le T}.
$$
We fix any $\nu\in\calv$ and
recall that $W$ is a Brownian motion under the probability
$\bar\P_\nu$. Taking the conditional expectation given $\bar\calf_t$
we obtain
$$
\bar Y_t'\ge
\bar \E_\nu\left[ \xi\,1_{\eta\ge T}+\int_{t\wedge \eta}^{T\wedge \eta}f_s\,ds +
h_\eta\,1_{t<\eta< T}\,\bigg|\, \bar\calf_t\right],
$$
where we have used again the fact that $\eta\ne T$ $\bar\P$-a.s. and hence
$\bar\P_\nu$-a.s.
Since $\nu$ was arbitrary in $\calv$
it follows that $
\bar Y_t'\ge
    J_t$ and in particular $
\bar Y_t\ge
    J_t$.

Next we prove the opposite inequality.
Comparing   \eqref{dualrepresentation} with
 \eqref{dualrepresentationpenalized}, since $\calv_n\subset\calv$
 it follows that $J_t^n\le J_t$. By the previous lemma
 we deduce that $\bar Y_t^n\le J_t$ and since
 $
   \bar Y^n\to
   \bar Y$ in $\cals^2$ we conclude that $\bar Y_t\le J_t$.
\qed

\noindent {\bf Conclusion of the proof of Theorem \ref{BSDEwellposed}.}
It remained to be shown that the solution $(\bar Y, \bar Z, \bar U, \bar K)$
constructed above is minimal.
Let
$(\bar Y', \bar Z', \bar U', \bar K')$
be any other solution to \eqref{BSDEconstrained}-\eqref{constraint}.
In the previous proof it was shown that, for every $t\in [0,T]$,
 $
\bar Y_t'\ge
    J_t$ $\bar\P$-a.s. Since we know from  Theorem \ref{dualproblem}
     that $
\bar Y_t=
    J_t$ we deduce that
    $
\bar Y_t'\ge \bar Y_t$. Since both processes are c\`adl\`ag,
this inequality holds for every $t$, up to a $\bar\P$-null set.
\qed



\begin{thebibliography}{11}

\bibitem{Bech} Becherer, D. Bounded solutions to backward SDE's with jumps for
utility optimization and indifference hedging.
Ann. Appl. Probab.  16  (2006),  no. 4, 2027-2054.


\bibitem{BeLi}
Bensoussan, A. and Lions, J.L.
Applications des in{\'e}quations variationnelles en contr{\^o}le stochastique.
  Dunod, 1978.


\bibitem{ElKh} Elie, R. and  Kharroubi I.
Adding constraints to BSDEs with jumps:
an alternative to multidimensional reflections.
ESAIM: Probability and Statistics   18 (2014),   233-250.

\bibitem{ElKh-b} Elie, R. and Kharroubi, I.
 BSDE representations for optimal switching problems
with controlled volatility.
Stoch. Dyn. 14  (2014),  1450003 (15 pages).




\bibitem{ElKetal97}
 El Karoui, N.; Kapoudjian, C.; Pardoux, E.; Peng, S.; Quenez, M. C. Reflected solutions of backward SDE's, and related obstacle problems for PDE's. Ann. Probab.  25  (1997),  no. 2, 702-737.

\bibitem{FuPh}  Fuhrman, M. and Pham, H.
Randomized and backward SDE
representation for optimal control of
non-markovian SDEs. Preprint
 arXiv:1310.6943, to appear on Ann. Appl.  Probab.

\bibitem{ja}
Jacod, J.   Multivariate point processes: predictable projection,
Radon-Nikodym derivatives, representation of martingales.
Z. Wahrscheinlichkeitstheorie und Verw. Gebiete   31 (1975), 235-253.


\bibitem{khaetal10} Kharroubi I., J. Ma, H. Pham, and J. Zhang.
Backward SDEs with constrained jumps and Quasi-variational inequalities.
  Ann. Probab.   38 (2010), 794-840.


\bibitem{khph}
Kharroubi, I. and  Pham H. (2012).  Feynman-Kac representation for Hamilton-Jacobi-Bellman IPDE.  Preprint  arXiv:1212.2000, to appear on   Ann.  Probab.

\end{thebibliography}
\end{document}